\documentclass[a4paper]{article}
\usepackage{amsmath,amssymb,amsfonts}
\usepackage{colordvi}
\usepackage[all]{xy}

\def\un{{\bf 1}}
\def\dom{\backslash}
\def\sur{\overline}

\def\mpn{\medskip\par\noindent}

\def\mmpn{\vskip 1em minus 1em\par\noindent}

\def\smp{\smallskip\par}

\newcommand{\flh}[2]{\mathop{\hbox to  4ex{\rightarrowfill}}_{#2}^{#1}
\limits}

\def\link{{\bf -}}

\def\Id{{\rm id}}
\def\Res{{\rm Res}}
\def\Ind{{\rm Ind}}
\def\Hom{{\rm Hom}}

\def\Aut{{\rm Aut}}

\def\Out{{\rm Out}}
\def\Inn{{\rm Inn}}
\def\Mod{{\rm Mod}}
\def\Inf{{\rm Inf}}
\def\Def{{\rm Def}}

\def\Iso{{\rm Iso}}
\def\Conj{{\rm Conj}}

\def\Ker{{\rm Ker}}

\def\Defres{{\rm Defres}}

\def\Indinf{{\rm Indinf}}

\def\Tr{{\rm Tr}}

\def\min{{\rm min}}

\def\ls#1#2{{\,^{#1}\!#2}}

\def\Z{\mathbb{Z}}

\def\pf{\par\bigskip\noindent{\bf Proof~: }}
\def\endpf{\nolinebreak~\leaders\hbox to 1em{\hss\
\hss}\hfill~\raisebox{.5ex}
{\framebox[1ex]{}}\par\bigskip}

\makeatletter
\renewenvironment{enumerate}{\ifnum \@enumdepth >3 \@toodeep\else
       \advance\@enumdepth \@ne
       \edef\@enumctr{enum\romannumeral\the\@enumdepth}\list
       {\csname  label\@enumctr\endcsname}{\setlength{\topsep}{1ex}
\setlength{\itemsep}{0 pt}\usecounter
         {\@enumctr}\def\makelabel##1{\hss\llap{##1}}}\fi}{\endlist}

\def\@seccntformat#1{\csname the#1\endcsname.\quad}

\def\section{\pagebreak[3]\setcounter{prop}{0}\setcounter{equation}{0}\@startsection{section}{1}{\z@}{4ex plus  6ex}{2ex}{\center\reset@font \large\bf}}

\makeatother

\def\theprop{\thesection.\arabic{prop}}

\newenvironment{enonce}[1]{\pagebreak[3]\refstepcounter{prop}\mmpn
{{\bf  \thesection.\arabic{prop}.\ #1.}}\begin{it} }{\end{it}\smp}
\def\thesection{\arabic{section}}

\newcommand{\result}[1]{\begin{enonce}{#1}}
\newcommand{\fresult}{\end{enonce}}

\begin{document}

\centerline{\Large\bf Vanishing evaluations of simple functors}
\vspace{.5cm}
\centerline{Serge Bouc, Radu Stancu, and Jacques Th\'evenaz}
\vspace{1cm}

\begin{footnotesize}
{\bf Abstract~:} The classification of simple biset functors is known, but the evaluation of a simple biset functor at a finite group $G$ may be zero. We investigate various situations where this happens, as well as cases where this does not occur. We also prove a closed formula for such an evaluation under some restrictive conditions on~$G$.

\bigskip\par

{\bf AMS Subject Classification~:} 19A22, 20C20.\par
{\bf Key words~:} biset, Burnside ring, simple module. 
\end{footnotesize}


\section{Introduction}
\noindent
Let $k$ be a field.
The {\it biset category\/} $k\cal C$ is the $k$-linear category whose objects are finite groups, with morphisms $\Hom_{k\cal C}(H,G)=kB(G,H)$, where $B(G,H)$ is the Burnside group of $(G,H)$-bisets and $kB(G,H)=k\otimes_{\Z}B(G,H)$.
A {\it biset functor\/} is a $k$-linear functor from $k\cal C$ to the category $k$-$\Mod$ of $k$-vector spaces. The category of biset functors is an abelian category and is used in various ways in representation theory, see \cite{Bo2}.\par

The classification of simple biset functors was obtained in~\cite{Bo1}. They are parametrized by equivalence classes of pairs $(H,V)$, where $H$ is a finite group and $V$ is a simple $k\Out(H)$-module. We write $S_{H,V}$ for the simple functor associated to the pair $(H,V)$.
However, the problem of describing the evaluation of simple functors at specific finite groups is much harder. In the present paper, we consider the question of the possible vanishing of such evaluations. We are interested in both questions of vanishing and non-vanishing.\par

This question is related to the problem of describing all simple modules for the double Burnside ring $kB(G,G)$, because any simple $kB(G,G)$-module has the form $S_{H,V}(G)$ for some $(H,V)$, and conversely any evaluation $S_{H,V}(G)$ is either zero or a simple $kB(G,G)$-module. We refer to \cite{BST} and \cite{BD} for this related question.\par

For other types of functors, there are explicit formulas for the evaluations of simple functors. This holds in particular for Mackey functors for a fixed finite group~$G$ (see Proposition~8.8 in~\cite{TW}) and for global Mackey functors and inflation functors (see Theorem~2.6 in~\cite{We}). Thus the vanishing of evaluations of such simple functors can be checked, at least in principle, by watching directly the formula. The situation is much more complicated for general biset functors, whenever both inflation and deflation are present (as well as restriction and induction). No closed formula for the evaluation of a simple functor is known. The most general known results are, on the one hand, a description of $S_{H,V}(G)$ as the image of a suitable linear map (see Theorem~4.3.20 in~\cite{Bo2}), and on the other hand, a formula for its dimension in terms of the rank of a suitable bilinear form (see Theorem~7.1 in~\cite{BST}). But neither of those results allows for an easy way to determine whether or not the evaluation is zero. The essential purpose of the present paper is to give some answers to this question.\par

In Section 3, we give some easy conditions which guarantee the non-vanishing $S_{H,V}(G)\neq 0$. Then we prove in Section~4 a general criterion, which has the disadvantage of being difficult to apply. In Section~5, we describe a suitable subquotient of $S_{H,V}(G)$, which implies a non-vanishing condition. Finally in Section~6, we prove that a closed formula for the evaluation at $G$ of a simple functor exists under some restrictive conditions on~$G$ and we immediately deduce a criterion for the vanishing of this evaluation. Various special cases can then be handled, as shown in Section~7.


\section{Preliminaries}
\noindent
We review some known facts about biset functors. For more details, we refer to~\cite{Bo1} and \cite{Bo2}. Given two finite groups $G$ and $H$, the Burnside group $B(G,H)$ is the Grothendieck group of the category of finite $(G,H)$-bisets and $kB(G,H)=k\otimes_{\Z}B(G,H)$. In particular, $kB(G,G)$ is a finite dimensional $k$-algebra, called the {\it double Burnside ring\/} of~$G$.\par

A {\it section\/} of a finite group $G$ is a pair $(T,S)$ of subgroups of $G$ such that $S$ is a normal subgroup of~$T$. In that case, the group $T/S$ is called a {\it subquotient\/} of~$G$. We write $H\sqsubseteq G$ when the group $H$ is isomorphic to a subquotient of~$G$ and we write $H\sqsubset G$ if $H\sqsubseteq G$ and $H\not\cong G$ (hence $|H|<|G|$).
We also write $N_G(T,S)$ for the normalizer of the section, that is, the set of all $g\in G$ such that $gTg^{-1}=T$ and $gSg^{-1}=S$.
If $(T,S)$ is a section of $G$, then there are elementary bisets $\Res_T^G$, $\Def_{T/S}^T$, $\Ind_T^G$, $\Inf_{T/S}^T$, and their composites $\Defres_{T/S}^G$ and $\Indinf_{T/S}^G$ (see Section~2.3 in~\cite{Bo2}). Also any group isomorphism $\sigma:G\to G'$ defines a $(G',G)$-biset $\Iso_\sigma$. Given finite groups $G$ and $H$, any transitive $(G,H)$-biset has the form $\Indinf_{B/A}^G\, \Iso_\sigma\,\Defres_{T/S}^H$, where $(B,A)$ is a section of $G$, $(T,S)$ is a section of $H$, and $\sigma:T/S\to B/A$ is a group isomorphism (see Lemma~3 in~\cite{Bo1} or Lemma~2.3.26 in~\cite{Bo2}).\par

 If $kI(G,G)$ is the ideal of $kB(G,G)$ generated by all $(G,G)$-bisets which factorize through a proper subquotient of~$G$, then $kB(G,G)/kI(G,G)\cong k\Out(G)$, where $\Out(G)=\Aut(G)/\Inn(G)$ is the group of outer automorphisms of~$G$. In particular, any $k\Out(G)$-module can be viewed as a $kB(G,G)$-module, with $kI(G,G)$ acting by zero.\par

The {\it biset category\/} $k\cal C$ is the $k$-linear category whose objects are finite groups, with morphisms $\Hom_{k\cal C}(H,G)=kB(G,H)$ (note that a $(G,H)$-biset is a morphism from $H$ to~$G$). The composition of morphisms, which we often write $\circ$, is the $k$-linear extension of the usual products of bisets $U\times_HV$. Recall that, if $U$ is a $(G,H)$-biset and $V$ is an $(H,L)$-biset, then $U\times_HV$ is a $(G,L)$-biset in the obvious way.\par

A {\it biset functor\/} is a $k$-linear functor from $k\cal C$ to the category $k$-$\Mod$ of $k$-vector spaces. The category of all such biset functors is abelian. A biset functor is called {\it simple\/} if it is non-zero and has no proper non-zero subfunctor. Recall the classification of simple functors
(see Section~4 in~\cite{Bo1} or Section~4.3 in~\cite{Bo2}).

\result{Proposition} \label{parametrization}
Let $S$ be a simple biset functor, let $H$ be a group of minimal order such that $S(H)\neq 0$, and let $V=S(H)$.
\begin{enumerate}
\item $H$ is unique up to isomorphism.
\item The ideal $kI(H,H)$ acts by zero on~$V$ and $V$ is a $k\Out(H)$-module.
\item $V$ is a simple $k\Out(H)$-module.
\item If $S(G)\neq0$ for some finite group $G$, then $H\sqsubseteq G$. 
\end{enumerate}
\fresult

This provides a {\it parametrization\/} of simple functors by (equivalence classes of) pairs $(H,V)$ where $H$ is a finite group and $V$ is a simple $k\Out(H)$-module. We write $S_{H,V}$ for the simple functor as in the statement, so that $H$ is its minimal group and $S_{H,V}(H)=V$.\par

We shall need a direct description of simple functors as quotients of suitable standard functors and we now recall this construction, which appears in~\cite{BST}. Let us fix a finite group $H$ and consider the representable functor $kB(-,H)$. For every finite group $G$, define
$$kI(G,H):=\sum_{J\sqsubset H}kB(G,J)B(J,H) \,.$$
Then $kI(-,H)$ is a subfunctor of $kB(-,H)$ and we define
$$k\sur B(-,H)=kB(-,H)/kI(-,H) \,.$$
For any finite group $G$, the evaluation $k\sur B(G,H)$ has a natural structure of right $k\Out(H)$-module, because the right action of $kI(H,H)$ is zero. This structure depends on~$G$ and we need to describe it more precisely.\par

We let $\Sigma_H(G)$ be the set of all sections $(T,S)$ of $G$ such that $T/S\cong H$, and we let $[\Sigma_H(G)/G]$ be a set of representatives of $G$-orbits in~$\Sigma_H(G)$. For every $(T,S)\in\Sigma_H(G)$, we choose an isomorphism $\sigma_{T,S}:H\to T/S$. The group $N_G(S,T)$ acts by conjugation on~$T/S$ and therefore $\sur N_G(S,T)=N_G(T,S)/T$ maps into the group $\Out(T/S)$. We use the isomorphism $\sigma_{T,S}$ to transport the image of this map to a subgroup of $\Out(H)$, that is, we define $\Gamma_G(T,S)$ to be the subgroup of $\Out(H)$ consisting of all elements induced by automorphisms $\sigma_{T,S}^{-1} \, \Conj_g \, \sigma_{T,S}$, where $g\in N_G(T,S)$. If $\sigma_{T,S}$ is replaced by $\sigma_{T,S}\, \alpha$ where $\alpha\in\Aut(H)$, then $\Gamma_G(T,S)$ is replaced by $\sur\alpha^{-1}\, \Gamma_G(T,S)\,\sur\alpha$, where $\sur\alpha\in\Out(H)$ is the class of~$\alpha$. Thus the conjugacy class of $\Gamma_G(T,S)$ only depends on~$(T,S)$.

\result{Lemma} \label{module-structure}
The $k$-space $k\sur B(G,H)$ is a permutation right $k\Out(H)$-module decomposing as follows~:
\begin{eqnarray*}
k\sur B(G,H) \;&=& {\displaystyle\bigoplus_{(T,S)\in[\Sigma_H(G)/G]} \sur{\Indinf_{T/S}^G}\circ k\sur B(T/S,H)} \\
&\cong& {\displaystyle\bigoplus_{(T,S)\in[\Sigma_H(G)/G]} k[\Gamma_G(T,S)\dom\Out(H)]} \,,
\end{eqnarray*}
where $\sur{\Indinf_{T/S}^G}$ denotes the image of~$\Indinf_{T/S}^G$ in~$k\sur B(G,H)$.
\fresult

\pf
By Lemma~7.2 in~\cite{BST}, $k\sur B(G,H)$ has a basis consisting of the set of elements of the form $\Indinf_{T/S}^G\circ\Iso_\sigma$, where $(T,S)\in[\Sigma_H(G)/G]$, and where $\sigma:H\to T/S$ runs over all isomorphisms, up to left composition by conjugation by elements of $N_G(T,S)$ (because $\Indinf_{T/S}^G=\Indinf_{T/S}^G\circ \Conj_g$ whenever $g\in N_G(T,S)$). This provides the first decomposition of the statement. Now for any fixed section $(T,S)$, we have a fixed isomorphism $\sigma_{T,S}:H\to T/S$ and we obtain a permutation right $k\Out(H)$-module
$$\sur{\Indinf_{T/S}^G}\circ k\sur B(T/S,H)
= \sur{\Indinf_{T/S}^G} \circ \sur{\Iso_{\sigma_{T,S}}} \circ k\Out(H) \,.$$
The generator $\sur{\Indinf_{T/S}^G} \circ \sur{\Iso_{\sigma_{T,S}}}$ has  $\Gamma_G(T,S)$ as a stabilizer in~$\Out(H)$. The result follows.
\endpf

For any left $k\Out(H)$-module~$V$, we define the functor
$$\sur L_{H,V}=k\sur B(-,H)\otimes_{k\Out(H)}V \,.$$
This has a subfunctor $\sur J_{H,V}$ defined as follows (see Remark~4.5 in~\cite{BST})~:
$$\sur J_{H,V}(X)=\Big\{\sum_i\sur\phi_i\otimes v_i \in \sur L_{H,V}(X) \mid \forall\psi\in kB(H,X), \sum_i (\sur{\psi\circ\phi}_i)\cdot v_i=0 \Big\} \,,$$
where $\phi_i\in kB(X,H)$ and $\sur\phi_i$ denotes its image in $k\sur B(X,H)$.
When $V$ is a simple $k\Out(H)$-module, we obtain the following result.

\result{Proposition} \label{L/J}
Suppose that $V$ is a simple $k\Out(H)$-module.
\begin{enumerate}
\item $\sur J_{H,V}$ is the unique maximal subfunctor of $\sur L_{H,V}$.
\item $\sur L_{H,V}/\sur J_{H,V}\cong S_{H,V}$.
\item For every finite group $G$ and for any fixed non-zero element $v\in V$,
$$\sur J_{H,V}(G)=\Big\{\sur\phi\otimes v \in \sur L_{H,V}(G) \mid \forall\psi\in kB(H,G), (\sur{\psi\circ\phi})\cdot v=0 \Big\} \,.$$
\item For every finite group $G$, $\sur L_{H,V}(G)$ is generated by all elements of the form $\sur{\Indinf_{T/S}^G\circ\Iso_{\sigma_{T/S}}}\otimes v$, where $v\in V$ and where $(T,S)$ is a section of $G$ such that $T/S\cong H$ and $\sigma_{T/S}:H\to S/T$ is a fixed isomorphism.
\end{enumerate}
\fresult

\pf
(1) and (2) follow from Proposition~4.4 in~\cite{BST}. Since $V$ is generated by any of its non-zero elements~$v$, (3) is a consequence of the description of $\sur J_{H,V}(X)$ given above. Finally (4) is a consequence of Lemma~\ref{module-structure} above.
\endpf


\section{Some easy cases}
\noindent
We have seen that $S_{H,V}(G)$ vanishes if $H$ is not isomorphic to a subquotient of~$G$.
Also, $S_{H,V}(H)=V\neq0$. This is of course the starting point in our investigation of vanishing or non-vanishing of evaluations. The following is another elementary result.

\result{Lemma} \label{quotient} Let $S_{H,V}$ be a simple biset functor. If $H$ is isomorphic to a quotient group $G/N$ of $G$, then $S_{H,V}(G)\neq0$.
\fresult

\pf
We have bisets $\Inf_{G/N}^G$ and $\Def_{G/N}^G$ which satisfy $\Def_{G/N}^G\circ\Inf_{G/N}^G=\Id_{G/N}$. Thus we have maps
$$S_{H,V}(G/N)  \flh{\Inf_{G/N}^G}{} S_{H,V}(G)  \flh{\Def_{G/N}^G}{} S_{H,V}(G/N)$$
whose composite is the identity. Since $S_{H,V}(G/N)\neq0$ (because $G/N\cong H$ by assumption), we must have $S_{H,V}(G)\neq0$.
\endpf

The lemma suffices to obtain the following result for the evaluation at an abelian group.

\result{Proposition} \label{abelian} Let $S_{H,V}$ be a simple biset functor.
If $H\sqsubseteq G/[G,G]$, then $S_{H,V}(G)\neq0$. In particular, if
$H\sqsubseteq G$ and if $G$ is abelian (hence $H$ too), then $S_{H,V}(G)\neq0$.
\fresult

\pf
In view of the structure theorem for finite abelian groups, any subquotient of the finite abelian group~$G/[G,G]$ is isomorphic to a quotient of~$G/[G,G]$, hence to a quotient of~$G$. Then the result follows from Lemma~\ref{quotient}.
\endpf

Our purpose is to generalize Lemma~\ref{quotient} and we need the following notions.
The set of all sections of~$G$ is partially ordered by the relation $\preceq$ defined as follows~:
$(V,U)\preceq (T,S)$ if and only if $V\leq T$ and the inclusion $\alpha:V\to T$ induces an isomorphism $\sur\alpha:V/U\to T/S$ (or in other words $VS=T$ and $V\cap S=U$).

Two sections $(B,A)$ and $(T,S)$ are said to be {\em linked\/} if $(B\cap T,A\cap S)\preceq (B,A)$ and $(B\cap T,A\cap S)\preceq (T,S)$ (see 4.3.11 in~\cite{Bo2} or Section~2 in~\cite{BT2}). In that case, the composition of the canonical isomorphisms
$$\phi_{B/A,T/S} : T/S \flh{\sim}{} (B\cap T)/(A\cap S) \flh{\sim}{} B/A$$
maps $xS$ to $xA$ for every $x\in B\cap T$ and is called the {\em isomorphism induced by the linking\/}. We write $(B,A)\link(T,S)$ whenever $(B,A)$ and $(T,S)$ are linked.

\result{Lemma} \label{linking} Let $(B,A)$ and $(T,S)$ be two sections of~$G$. The following are equivalent~:
\begin{enumerate}
\item[a)] $(B,A)\link(T,S)$,
\item[b)] $|B/A| = |T/S|$, $S(B\cap T)=T$, and $S(A\cap T)=S$.
\end{enumerate}
\fresult

Of course, the last equality is equivalent to $A\cap T\leq S$, but we shall need below the equality as stated.

\pf
We know that $(B,A)$ and $(T,S)$ generate a butterfly diagram, as in Lemma~2.3 of~\cite{BT2},
and the two sections
$$(A(B\cap S),A(B\cap T))\quad\text{ and } \quad (S(B\cap T),S(A\cap T))$$
are linked.
If now $S(B\cap T)=T$ and $S(A\cap T)=S$, then the second section is $(T,S)$.
Thus $(T,S)$ is linked to $(A(B\cap S),A(B\cap T))$, which is a section of~$B/A$.
If moreover, $|B/A| = |T/S|$, then this section of $B/A$ cannot be proper and must be~$(B,A)$.
It follows that $(T,S)$ is linked to $(B,A)$ (i.e. the butterfly diagram collapses to a linking).

If conversely $(B,A)\link(T,S)$, then $B/A\cong T/S$, hence $|B/A| = |T/S|$.
Moreover, the linking implies that $S(B\cap T)=T$, and also that $A\cap T=A\cap S$, so that $S(A\cap T)=S$.
\endpf

\result{Proposition} \label{not-linked} Let $S_{H,V}$ be a simple biset functor.
Suppose that $H$ is isomorphic to $T/S$, where $(T,S)$ is a section of $G$ such that, for every $g\in G$ with $g\notin T$, the conjugate section $(\ls gT,\ls gS)$ is not linked to~$(T,S)$.
Then $S_{H,V}(G)\neq0$.
\fresult

It is clear that if $S$ is a normal subgroup $S$ of $G$, then the section $(G,S)$ satisfies the assumption (because in that case there is no $g\notin G$), so Proposition~\ref{not-linked} actually generalizes Lemma~\ref{quotient}.

\pf
As in the proof of Lemma~\ref{quotient}, we consider the maps
$$S_{H,V}(T/S)  \flh{\Indinf_{T/S}^G}{} S_{H,V}(G) 
\flh{\Defres_{T/S}^G}{} S_{H,V}(T/S)$$
and we want to prove that the composite is the identity.
This will then force $S_{H,V}(G)$ to be non-zero since
$S_{H,V}(T/S) \cong S_{H,V}(H)\neq0$.

By the generalized Mackey formula (see Proposition~A.1 in~\cite{BT1}),
the composite above decomposes as a sum indexed by double cosets representatives $g\in[T\dom G/T]$.
There is one double coset, indexed by an element of~$T$ which can be chosen to be~$1_G$,
and the corresponding term is the identity.
We show that all the other terms vanish. Any such term has the form
$$\Indinf_{X/Y}^{T/S} \;\Iso_\phi \; \Conj_g \;
\Defres_{S(\ls gT\cap T)/S(\ls gS\cap T)}^{T/S}
$$
for some subquotient $X/Y$ of $T/S$ and some group isomorphism~$\phi$.
Since $(\ls gT,\ls gS)$ is not linked to~$(T,S)$ by assumption,
and since $|\ls gT/\ls gS| = |T/S|$, Lemma~\ref{linking} tells us that the section
$(S(\ls gT\cap T),S(\ls gS\cap T))$ must be a proper section of~$T/S$.
Therefore the group
$S(\ls gT\cap T)/S(\ls gS\cap T)$ has order strictly smaller than the order of $T/S\cong H$. The functor $S_{H,V}$ vanishes on such a group and so the term above factors through zero.
\endpf

The special case where $S=1$ is worth mentioning.

\result{Corollary} \label{self-normalizing} 
Suppose that $H$ is isomorphic to a subgroup $T$ of~$G$ such that $T=N_G(T)$. Then $S_{H,V}(G)\neq0$.
\fresult

We now show that the assumption of Proposition~\ref{not-linked} holds in particular for a section $(N_G(S),S)$ where $S$ is an expansive subgroup of~$G$. Recall that a subgroup $S$ of~$G$ is called {\it expansive\/} in~$G$ if, for every $g\notin N_G(S)$, the subgroup $S(\ls gS\cap N_G(S))/S$ has a non-trivial core in the group $N_G(S)/S$, in other words, there exists a normal subgroup $M$ of~$N_G(S)$ contained in~$S(\ls gS\cap N_G(S))$ and containing~$S$ properly. This notion is defined and used in~\cite{Bo2} and in~\cite{BT2}. In particular,
the subgroup $S(\ls gS\cap N_G(S))$ contains $S$ properly, and this implies, by Lemma~\ref{linking}, that the section $(\ls gN_G(S),\ls gS)$ is not linked to~$(N_G(S),S)$, so that Proposition~\ref{not-linked} applies. This proves the following corollary.

\result{Corollary} \label{expansive} Let $S_{H,V}$ be a simple biset functor.
If $H$ is isomorphic to $N_G(S)/S$, where $S$ is an expansive subgroup of~$G$, then $S_{H,V}(G)\neq0$.
\fresult

It is clear that any normal subgroup $S$ of $G$ is expansive in~$G$ (because in that case there is no $g\notin N_G(S)$), so again Lemma~\ref{quotient} is a special case of Corollary~\ref{expansive}.

For example, the Mathieu group $G=M_{11}$ has a subgroup~$S$, isomorphic to~$A_6$, which is expansive in~$G$ and such that $N_G(S)/S$ has order~2 (in fact $N_G(S)=M_{10}$). It follows that $S_{C_2,k}(M_{11})\neq 0$,
independently of the characteristic of~$k$.


\section{A general criterion}
\noindent
Let $S_{H,V}$ be a simple biset functor and let $G$ be a finite group. The analysis of the evaluation $S_{H,V}(G)$ involves the set $\Sigma_H(G)$ of all sections $(T,S)$ of $G$ such that $T/S\cong H$, because $S_{H,V}$ is a quotient of $\sur L_{H,V}$ and the evaluation $\sur L_{H,V}(G)$ involves those sections (see Proposition~\ref{L/J}). We may assume that $H\sqsubseteq G$, that is, $\Sigma_H(G)\neq\emptyset$.\par

Now we come to a criterion for the vanishing of the evaluation of a simple functor. It gives a general answer to our main question, although it is rather hard to use it in practice. A similar result appears in Theorem~7.1 of~\cite{BD}.

\result{Theorem} \label{criterion}
Let $S_{H,V}$ be a simple biset functor and let $G$ be a finite group. For every $(T,S)\in\Sigma_H(G)$, fix an isomorphism $\sigma_{T/S}:H\to T/S$. The following are equivalent~:
\begin{enumerate}
\item $S_{H,V}(G)=0$.
\item For any $(B,A),(T,S)\in\Sigma_H(G)$,
the action on~$V$ of the automorphism
$$\sum_{{\scriptstyle g\in[B\dom G/T]}\atop{\scriptstyle (B,A)\link\ls g(T,S)}}
\sigma_{B/A}^{-1} \; \phi_{B/A,\ls gT/\ls gS} \;  \Conj_g \; \sigma_{T/S}$$
is zero, where $\phi_{B/A,\ls gT/\ls gS}:\ls gT/\ls gS \to B/A$ denotes the isomorphism induced by the linking $(B,A)\link\ls g(T,S)$.
\end{enumerate}
\fresult

\pf
By Proposition~\ref{L/J}, we have $S_{H,V}(G)=0$ if and only if $\sur L_{H,V}(G)=\sur J_{H,V}(G)$. By Proposition~\ref{L/J} again, the latter equality holds if and only if $(\sur{\psi\circ\phi})$ acts by zero on~$V$ for all $\phi\in kB(G,H)$ and  all $\psi\in kB(H,G)$. Since we have passed to the quotient by all morphisms factorizing below~$H$, we can assume that
$$\phi = \Indinf_{T/S}^G\circ\Iso_{\sigma_{T/S}} \quad \text{and}\quad \psi = \Iso_{\sigma_{B/A}}^{-1} \circ\Defres_{B/A}^G$$
where $(B,A),(T,S)\in\Sigma_H(G)$. Then the generalized Mackey formula applies to~$\Defres_{B/A}^G\circ\Indinf_{T/S}^G$ (see Proposition~A.1 in~\cite{BT1}), indexed by double cosets representatives $g\in[B\dom G/T]$. But we are interested in the image $\sur{\psi\circ\phi}$ in
$$kB(H,H)/kI(H,H)\cong k\Out(H)$$
and all terms in the formula factorize through a group isomorphic to a proper subquotient of~$H$, except those indexed by an element $g\in[B\dom G/T]$ such that $(B,A)\link\ls g(T,S)$, where $\ls g(T,S)=(\ls gT, \ls gS)$ denotes the conjugate section. For such an element $g$, we are left with the $(B/A,T/S)$-biset $\Iso_{\phi_{B/A,\ls gT/\ls gS}} \circ \Iso_{\Conj_g}$. Composing with $\sigma_{T/S}$ and $\sigma_{B/A}^{-1}$, we see that the action on~$V$ of $(\sur{\psi\circ\phi})$ is equal to the action of the automorphism
$$\sum_{{\scriptstyle g\in[B\dom G/T]}\atop{\scriptstyle (B,A)\link\ls g(T,S)}} \sigma_{B/A}^{-1} \,
\phi_{B/A,\ls gT/\ls gS} \,  \Conj_g\,\sigma_{T/S} \,.$$
Thus the condition is that this sum must act by zero on~$V$.
\endpf

Note that if $(B,A)\link\ls g(T,S)$, then the isomorphism $\phi_{B/A,\ls gT/\ls gS} \, \Conj_g$ induced by the linking is given by the $(B/A,T/S)$-biset $A\dom BgT/S$. This appears explicitly in the proof of the generalized Mackey formula in~\cite{BT1}, but we do not need this here.\par


\section{Minimal sections}
\noindent
Given finite groups $H$ and $G$, we let again $\Sigma_H(G)$ be the set of all sections $(T,S)$ of $G$ such that $T/S\cong H$.
A section $(T,S)\in\Sigma_H(G)$ will be called {\em minimal\/} if it is minimal with respect to the partial order $\preceq$ defined in Section~3. In that case, if $(T,S)$ is linked to $(B,A)$, then ${(B\cap T,A\cap S)} = (T,S)$, that is, $T\leq B$ and $S\leq A$ (and also $TA=B$ and $T\cap A=S$ because of the linking). We write $\Sigma_H(G)^\min$ for the subset of $\Sigma_H(G)$ consisting of minimal sections. Clearly $G$ acts by conjugation on $\Sigma_H(G)$ and $\Sigma_H(G)^\min$ and we let $[\Sigma_H(G)^\min/G]$ denote a set of representatives of $G$-orbits in~$\Sigma_H(G)^\min$.

\result{Lemma} \label{Frattini}
Let $(T,S)\in\Sigma_H(G)$, let $f:T\to H$ be a surjective group homomorphism with kernel~$S$, and let $\Phi(T)$ be the Frattini subgroup of~$T$ (that is, the intersection of all maximal subgroups of~$T$). The following are equivalent.
\begin{enumerate}
\item $(T,S)$ is minimal.
\item $S\leq\Phi(T)$.
\item $f$ induces an isomorphism $T/\Phi(T)\flh{\sim}{} H/\Phi(H)$.
\end{enumerate}
\fresult

\pf
If $H=\un$, then the only minimal section in $\Sigma_\un(G)$ is $(\un,\un)$ and the result follows easily. Assume that $H\neq\un$, that is $S<T$.
Suppose $(T,S)$ is not minimal and let $(B,A)\in\Sigma_H(G)$ such that $(B,A)\prec(T,S)$. Then $B<T$ and there is some maximal subgroup $M$ of $T$ containing~$B$. It follows that $T=BS=MS$, so $S\not\leq\Phi(T)$. Conversely, if $S\not\leq\Phi(T)$, there is some maximal subgroup $M$ of $T$ which does not contain~$S$. Then $MS=T$, so $(M,M\cap S)\prec(T,S)$ and $(T,S)$ is not minimal. The proof of the equivalence of (2) and~(3) is easy and is left to the reader.
\endpf

Recall that, for every section $(T,S)\in\Sigma_H(G)$, we have set $\sur N_G(S,T)=N_G(S,T)/T$ and we have fixed an isomorphism $\sigma_{T/S}:H\to T/S$. This allows us to view any $k\Out(H)$-module $V$ as a $k[\sur N_G(T,S)]$-module, as follows~:
$$\sur g\cdot v= \sur{\sigma_{T/S}^{-1} \, \Conj_g \, \sigma_{T/S}}\cdot v \,,\quad g\in N_G(T,S)\,,v\in V \,,$$
where the bar denotes the class in~$\Out(H)$ of the automorphism in~$\Aut(H)$.
By Proposition~\ref{L/J}, we know that $\sur L_{H,V}(G)$ is generated as a $k$-vector space by all the elements of the form
$$\sur{\Indinf_{T/S}^G\circ\Iso_{\sigma_{T/S}}}\otimes v \,,\quad(T,S)\in\Sigma_H(G), \; v\in V \,.$$
Let $\sur L_{H,V}(G)^\min$ be the subspace of $\sur L_{H,V}(G)$ generated by all the elements of the form
$$\sur{\Indinf_{T/S}^G\circ\Iso_{\sigma_{T/S}}}\otimes v \,,\quad(T,S)\in\Sigma_H(G)^\min, \; v\in V \,.$$
Let also $S_{H,V}(G)^\min$ be the image of $\sur L_{H,V}(G)^\min$ under the canonical surjection
$$\pi:\sur L_{H,V}(G)\flh{}{} \sur L_{H,V}(G)/\sur J_{H,V}(G) \cong S_{H,V}(G) \,.$$
For any finite group $X$ and any $kX$-module $W$, we let $\Tr_1^X:W\to W$ be the $k$-linear map defined by $\Tr_1^X(w)=\sum_{x\in X} x{\cdot}w$ (relative trace), and we let $\Tr_1^X(W)$ denote its image.

\result{Theorem} \label{minimal-sections}
Let $S_{H,V}$ be a simple biset functor and let $G$ be a finite group.
With the notation above, there is a surjective $k$-linear map
$$\sur\tau: S_{H,V}(G)^\min \flh{}{} \bigoplus_{(T,S)\in[\Sigma_H(G)^\min/G]} \Tr_1^{\sur N_G(T,S)}(V) \,.$$
Hence the right hand side is isomorphic to a subquotient of $S_{H,V}(G)$.
\fresult

\pf
We have $\sur L_{H,V}(G)=k\sur B(G,H)\otimes_{k\Out(H)}V$ and we know how $k\sur B(G,H)$ decomposes, by Lemma~\ref{module-structure}. By tensoring with~$V$ this decomposition, we obtain
\begin{eqnarray*}
\sur L_{H,V}(G) \;=&{\displaystyle\bigoplus_{(T,S)\in[\Sigma_H(G)/G]}}&
\sur{\Indinf_{T/S}^G}\circ k\sur B(T/S,H) \otimes_{k\Out(H)}V \\
\cong& {\displaystyle\bigoplus_{(T,S)\in[\Sigma_H(G)/G]}}&
k[\Gamma_G(T,S)\dom\Out(H)] \otimes_{k\Out(H)}V \\
\cong& {\displaystyle\bigoplus_{(T,S)\in[\Sigma_H(G)/G]}}&
V_{\Gamma_G(T,S)}  \,, 
\end{eqnarray*}
where $V_{\Gamma_G(T,S)}$ denotes the $k$-space of coinvariants for the group $\Gamma_G(T,S)$ (i.e. the quotient of~$V$ by all elements of the form $(\gamma{-}1)v$, $\gamma\in\Gamma_G(T,S)$, $v\in V$).
By the proof of Lemma~\ref{module-structure}, we see that a generator
$\sur{\Indinf_{T/S}^G\circ\Iso_{\sigma_{T/S}}}\otimes v$ of $\sur L_{H,V}(G)$ is mapped to $1\otimes v \in k[\Gamma_G(T,S)\dom\Out(H)] \otimes_{k\Out(H)}V$, hence to the class of $v$ in $V_{\Gamma_G(T,S)}$.

By definition, $\Gamma_G(T,S)$ is the image of $\sur N_G(T,S)$ in~$\Out(H)$ (using the isomorphism $\sigma_{T/S}:H\to T/S$). Therefore $V_{\Gamma_G(T,S)}=V_{\sur N_G(T,S)}$ and we obtain
$$\sur L_{H,V}(G) = \bigoplus_{(T,S)\in[\Sigma_H(G)/G]} V_{\sur N_G(T,S)} \,.$$
Restricting to minimal sections, we obtain
$$\sur L_{H,V}(G)^\min \cong \bigoplus_{(T,S)\in[\Sigma_H(G)^\min/G]} V_{\sur N_G(T,S)}\,.$$
Now the relative trace map $\Tr_1^{\sur N_G(T,S)}$ induces a surjective $k$-linear map
$$\Tr_1^{\sur N_G(T,S)} : V_{\sur N_G(T,S)} \flh{}{} \Tr_1^{\sur N_G(T,S)}(V) \,,$$
and the direct sum of these maps defines a surjective $k$-linear map
$$\tau: \sur L_{H,V}(G)^\min \; \flh{}{}  \bigoplus_{(T,S)\in[\Sigma_H(G)^\min/G]} \Tr_1^{\sur N_G(T,S)}(V) \,,$$
mapping a generator $\sur{\Indinf_{T/S}^G\circ\Iso_{\sigma_{T/S}}}\otimes v$
to $\Tr_1^{\sur N_G(T,S)}(v)$.

Recall the canonical surjection
$$\pi:\sur L_{H,V}(G)\flh{}{} \sur L_{H,V}(G)/\sur J_{H,V}(G) \cong S_{H,V}(G) \,.$$
We claim that $\sur J_{H,V}(G) \cap \sur L_{H,V}(G)^\min \subseteq \Ker(\tau)$. It will follow that $\tau$ induces a surjective $k$-linear map
$$\sur\tau: S_{H,V}(G)^\min \flh{}{} \bigoplus_{(T,S)\in[\Sigma_H(G)^\min/G]} \Tr_1^{\sur N_G(T,S)}(V) \,,$$
proving the theorem.

In order to prove the claim, we let $x\in \sur J_{H,V}(G) \cap \sur L_{H,V}(G)^\min$ and we write
$$x=\sum_{(T,S)\in[\Sigma_H(G)^\min/G]} (\sur{\Indinf_{T/S}^G\circ\Iso_{\sigma_{T/S}}}) \otimes v_{T,S} \,,$$
where  $v_{T,S}\in V$ for every $(T,S)$.
By the description of $J_{H,V}(G)$ in Proposition~\ref{L/J}, we have
$$\sum_{(T,S)\in[\Sigma_H(G)^\min/G]} (\sur{\psi}\circ\sur{\Indinf_{T/S}^G\circ\Iso_{\sigma_{T/S}}})\cdot v_{T,S} =0$$
for all $\psi\in kB(H,G)$. Fix $(B,A)\in[\Sigma_H(G)^\min/G]$ and choose
$$\psi = \Iso_{\sigma_{B/A}}^{-1} \circ\Defres_{B/A}^G \,.$$
As in the proof of Theorem~\ref{criterion}, $\Defres_{B/A}^G\circ\Indinf_{T/S}^G$ decomposes according to the generalized Mackey formula (see Proposition~A.1 in~\cite{BT1}), indexed by double cosets representatives $g\in[B\dom G/T]$, and all terms in the formula factorize through a group isomorphic to a proper subquotient of~$H$, except those indexed by an element $g\in[B\dom G/T]$ such that $(B,A)\link\ls g(T,S)$. Since both $(B,A)$ and $\ls g(T,S)$ are minimal, the only possible linking is the identity, hence $(B,A)=\ls g(T,S)$. But we have chosen orbit representatives in~$\Sigma_H(G)^\min$, so $(B,A)=(T,S)$. Thus the sum over $(T,S)$ reduces to a single term, indexed by~$(B,A)$. Moreover, in the Mackey formula, we are left with the sum over all $g\in[B\dom G/B]$ such that $\ls g(B,A)=(B,A)$, that is, $g\in [N_G(B,A)/B]$. Therefore we obtain
$$\sum_{g\in [N_G(B,A)/B]}  (\Iso_{\sigma_{B/A}}^{-1} \circ \Iso_{\Conj_g} \circ \Iso_{\sigma_{B/A}})
\cdot v_{B,A} =0 \,.$$
This is the action of the automorphism $\sigma_{B/A}^{-1} \, \Conj_g \, \sigma_{B/A}$, so by definition of the action of $\sur N_G(B,A)$, we obtain
$$\Tr_1^{\sur N_G(B,A)} (v_{B,A})=0 \,.$$
This holds for all $(B,A)\in \Sigma_H(G)^\min$ and therefore
\begin{eqnarray*}
\tau(x)&= \tau\bigl(\displaystyle\sum_{(T,S)\in[\Sigma_H(G)^\min/G]}
\sur{\Indinf_{T/S}^G\circ\Iso_{\sigma_{T/S}}}\otimes v_{T,S}\bigr) \\
&=\displaystyle\sum_{(T,S)\in[\Sigma_H(G)^\min/G]} \Tr_1^{\sur N_G(T,S)} (v_{T,S})=0 \,.
\end{eqnarray*}
Thus $x\in\Ker(\tau)$, proving the claim.
\endpf

Theorem~\ref{minimal-sections} immediately implies the following result about vanishing evaluations.

\result{Corollary} \label{minimal-vanishing}
Let $S_{H,V}$ be a simple biset functor and let $G$ be a finite group. If $\Tr_1^{\sur N_G(T,S)}(V) \neq0$ for some $(T,S)\in\Sigma_H(G)^\min$, then $S_{H,V}(G)\neq0$.
\fresult

This can be applied for instance in the following situation.

\result{Corollary} \label{trivial-Gamma}
Suppose that there exists a minimal section $(T,S)\in\Sigma_G(H)^\min$ such that $N_G(T,S)$ acts by inner automorphisms on~$T/S$. Suppose also that $|\sur N_G(T,S)|\neq0$ in~$k$. Then $S_{H,V}(G)\neq 0$.
\fresult

\pf
By assumption, the image of $\sur N_G(T,S)$ in~$\Out(T/S)$ is trivial.
Therefore $\Tr_1^{\sur N_G(T,S)}(V) = |\sur N_G(T,S)|\cdot V\neq0$ and the result follows from Corollary~\ref{minimal-vanishing}.
\endpf

Note that the first assumption holds in particular if $N_G(T,S)$ is equal to the centralizer $C_G(T/S)$ of $T/S$, because $N_G(T,S)$ acts trivially on~$T/S$ in this case. This applies in particular if $N_G(T,S)$ is abelian, improving the first statement of Proposition~\ref{abelian} in the case where $k$ has characteristic not dividing~$|G|$.


\section{A closed formula for some evaluations}
\noindent
As already mentioned in the introduction, there is no known closed formula for the evaluation $S_{H,V}(G)$ of a simple biset functor $S_{H,V}$. However, with suitable assumptions, such a formula exists.\par

For instance, if the $k\Out(H)$-module $V$ is {\em primitive\/}, in the sense defined on page~721 of~\cite{Bo1}, and if $k$ has characteristic zero, then
$$S_{H,V}(G) \cong \bigoplus_{(T,S)} V^{\sur N_G(T,S)} \,,$$
where $(T,S)$ runs over a suitable subset of $\Sigma_H(G)$. We refer to Proposition~20 in~\cite{Bo1} for more details. Of course, a criterion for the vanishing of $S_{H,V}(G)$ is immediately deduced in this case.\par

The purpose of this section is to prove a closed formula for the evaluation $S_{H,V}(G)$ under a suitable assumption on the structure of~$G$, more precisely when $\Sigma_H(G)^\min=\Sigma_H(G)$. Then this can be used to give a criterion for the vanishing of $S_{H,V}(G)$.

\result{Theorem} \label{closed-formula}
Let $S_{H,V}$ be a simple biset functor and let $G$ be a finite group. If every section $(T,S)\in\Sigma_H(G)$ is minimal, the map $\sur\tau$ of Theorem~\ref{minimal-sections} is an isomorphism and
$$S_{H,V}(G) \cong \bigoplus_{(T,S)\in[\Sigma_H(G)/G]} \Tr_1^{\sur N_G(T,S)}(V) \,.$$
\fresult

\pf
By Theorem~\ref{minimal-sections}, we already know that the map $\sur\tau$ is surjective, so we need to prove that it is injective. In other words, we have to show that
$$\Ker(\tau)\subseteq \sur J_{H,V}(G) \cap \sur L_{H,V}(G)^\min \,,$$
where $\tau$ is the map defined in the proof of Theorem~\ref{minimal-sections}, namely
$$\tau: \sur L_{H,V}(G)^\min \; \flh{}{}  \bigoplus_{(T,S)\in[\Sigma_H(G)^\min/G]} \Tr_1^{\sur N_G(T,S)}(V) \,.$$
Let $x\in \Ker(\tau)$ and write
$$x\;=\sum_{(T,S)\in[\Sigma_H(G)^\min/G]} \sur{\Indinf_{T/S}^G\circ\Iso_{\sigma_{T/S}}}\otimes v_{T,S}  \,,$$
where  $v_{T,S}\in V$ for every $(T,S)$.
Since $\tau(x)=0$, we have $\Tr_1^{\sur N_G(T,S)} (v_{T,S})=0$ for every $(T,S)\in\Sigma_H(G)^\min$.  Fix $(B,A)\in[\Sigma_H(G)^\min/G]$ and let
$$\psi_{B,A} = \Iso_{\sigma_{B/A}}^{-1} \circ\Defres_{B/A}^G \,.$$
Exactly the same computation as in the proof of Theorem~\ref{minimal-sections} above shows that
\begin{eqnarray*}
\sur{\psi_{B,A}}\cdot x&=&
\sum_{(T,S)\in[\Sigma_H(G)^\min/G]} (\sur{\psi_{B,A}}\circ\sur{\Indinf_{T/S}^G\,\Iso_{\sigma_{T/S}}})
\cdot v_{T,S} \\
&=&\Tr_1^{\sur N_G(B,A)} (v_{B,A}) =0 \,.
\end{eqnarray*}
But now, since $\Sigma_H(G)^\min=\Sigma_H(G)$ by assumption, $k\sur B(H,G)$ is generated by elements of the form $\sur{\Iso_\alpha\circ\psi_{B,A}}$ where $\psi_{B,A}$ is as above and $\alpha\in\Aut(H)$. Therefore
$\sur{\psi}\cdot x=0$ for all $\psi\in kB(H,G)$. But this means that $x\in \sur J_{H,V}(G)$, as was to be shown.
\endpf

\result{Corollary} \label{criterion-vanishing}
Let $S_{H,V}$ be a simple biset functor and let $G$ be a finite group.
Assume that $\Sigma_H(G)^\min=\Sigma_H(G)$. Then $S_{H,V}(G)=0$ if and only if $\Tr_1^{\sur N_G(T,S)}(V) =0$ for every $(T,S)\in\Sigma_H(G)$.
\fresult

Our next result gives a first application of Theorem~\ref{closed-formula}.

\result{Proposition} \label{p-groups}
Suppose that $G$ and $H$ are $p$-groups with the same sectional rank. Then
$$S_{H,V}(G) \cong \bigoplus_{(T,S)\in[\Sigma_H(G)/G]} \Tr_1^{\sur N_G(T,S)}(V) \,.$$
In particular $S_{H,V}(G)=0$ if and only if the action of $\sum_{g\in [N_G(T,S)/T]} \Conj_g$ on~$V$ is zero, for every $(T,S)\in\Sigma_H(G)$.
\fresult

\pf
Let $(T,S)\in\Sigma_H(G)$ and let $f:T\to H$ be a surjective group homomorphism with kernel~$S$.
Let $r$ be the sectional rank of~$H$ and let $1\leq H_2\triangleleft H_1\leq H$ such that $H_1/H_2$ is elementary abelian of rank~$r$.
Let $T_i=f^{-1}(H_i)$. Then $T_1/T_2$ is elementary abelian of rank~$r$, and this must be the largest possible rank of an elementary abelian quotient of~$T_1$, because the sectional rank of~$G$ is also~$r$. It follows that $\Phi(T_1)=T_2$. Since $\Phi(T_1)\leq \Phi(T)$ (because $T_1\leq T$ and $T$ is a $p$-group), we deduce that $S\leq T_2 =\Phi(T_1)\leq \Phi(T)$. By Lemma~\ref{Frattini}, this proves that the section $(T,S)$ is minimal. Thus Theorem~\ref{closed-formula} applies and yields the result.
\endpf


\section{The case of a single section}
\noindent
Theorem~\ref{closed-formula} can be applied to various cases to obtain a closed formula for the evaluation $S_{H,V}(G)$, hence a criterion for its vanishing. We concentrate here on cases where the set $\Sigma_H(G)$ reduces to a single conjugacy class, so that clearly every section $(T,S)\in\Sigma_H(G)$ is minimal. As before, $S_{H,V}$ denotes a simple biset functor and $G$ a finite group. We can assume that $H\sqsubseteq G$, since otherwise $S_{H,V}(G)=0$.

\result{Proposition} \label{only-section}
Suppose that $\Sigma_H(G)$ contains a unique section $(T,S)$ up to conjugation. Then
$$S_{H,V}(G) \cong \Tr_1^{\sur N_G(T,S)}(V) \,.$$
In particular $S_{H,V}(G)=0$ if and only if the action of $\sum_{g\in [N_G(T,S)/T]} \Conj_g$ on~$V$ is zero.
\fresult

\pf
This is a special case of Theorem~\ref{closed-formula}.
\endpf

There are many instances where $G$ has a subgroup $H$ such that $(H,\un)$ is the only section in $\Sigma_H(G)$ up to conjugation. Here are a few such cases.

\result{Corollary} \label{Hall}
Suppose that $H$ is a normal Hall subgroup of $G$. By the Schur-Zassenhaus theorem, we know that $G=H\rtimes Y$ for some subgroup~$Y$. Assume that $Y$ acts faithfully on~$H$. Then
$$S_{H,V}(G) \cong \Tr_1^Y (V) \,.$$
In particular $S_{H,V}(G)=0$ if and only if the action of $\sum_{y\in Y}\Conj_y$ on~$V$ is zero.
\fresult

\pf
Since $(|H|,|Y|)=1$ and $Y$ acts faithfully on~$H$, it is easy to see that $H$ cannot normalize a non-trivial subgroup of~$Y$. Therefore there is no section in $\Sigma_H(G)$ apart from $(H,\un)$ and Proposition~\ref{only-section} applies. Moreover, $N_G(H)/H=G/H\cong Y$.
\endpf

\result{Corollary} \label{Sylow} 
Suppose that $H$ is a Sylow $p$-subgroup of $G$ and that $H$ does not normalize any non-trivial $q$-subgroup, for every prime $q\neq p$.  Then
$$S_{H,V}(G) \cong \Tr_1^{N_G(H)/H}(V) \,.$$
In particular $S_{H,V}(G)=0$ if and only if the action of $\sum_{g\in [N_G(H)/H]} \Conj_g$ on~$V$ is zero.
\fresult

\pf
Let $(T,S)$ be a section of $G$ such that $T/S\cong H$. Then a Sylow $p$-subgroup of $T$ is isomorphic to $H$, hence conjugate to~$H$, and we may assume that it is equal to~$H$. Then $H$ normalizes $S$, hence also a Sylow $q$-subgroup $Q$ of~$S$ because $(|H|,|S|)=1$. Therefore $Q=1$ by assumption. Since this holds for every prime divisor $q$ of~$|S|$, it follows that $S=1$ and that $T=H$. Thus $\Sigma_H(G)$ reduces to the conjugacy class of $(H,\un)$ and Proposition~\ref{only-section} applies.
\endpf

\result{Corollary} \label{simple} 
Let $H$ be a non-abelian simple group and let $G=\Aut(H)$. Suppose that $V$ is a non-trivial $k\Out(H)$-module (so in particular $G/H\cong \Out(H)$ is a non-trivial group). Then $S_{H,V}(G)=0$.
\fresult

\pf
Let $(T,S)$ be a section of $G$ such that $T/S\cong H$. Then $HT/HS$ is isomorphic to a quotient of $T/S\cong H$ and to a subquotient of $G/H\cong \Out(H)$. Since $\Out(H)$ is solvable (the Schreier conjecture) and $H$ is non-abelian simple, we must have $HT/HS=\un$, hence $HT=HS$. It follows that $T=(H\cap T)S$, hence $T/S\cong(H\cap T)/(H\cap S)$. But $T/S\cong H$, while $(H\cap T)/(H\cap S)$ is a subquotient of~$H$. Therefore $H\cap T=H$ and $H\cap S=\un$. But then $H$ and $S$ are normal subgroups of~$T$ with $H\cap S=\un$, so $[S,H]=\un$, that is, $S$ acts trivially on~$H$. Since $S$ is isomorphic to a subgroup of $\Out(H)$, we must have $S=\un$. This proves that $(H,\un)$ is the only section in~$\Sigma_H(G)$ and therefore
$$S_{H,V}(G) \cong \Tr_1^{G/H}(V) = \Tr_1^{\Out(H)}(V) \subseteq V^{\Out(H)} =0 \,,$$
because the non-trivial simple module $V$ has no non-zero fixed point.
\endpf

If $k$ has characteristic zero, then under the assumptions of Corollary~\ref{simple}, the module $V$ is in fact primitive, in the sense mentioned at the beginning of Section~6, and this provides another approach of the result.

\result{Corollary} \label{central} 
Let $H$ be a non-abelian simple group and let $G$ be its universal central extension. Then $S_{H,V}(G)=V$ (hence non-zero).
\fresult

\pf
Let $Z$ be the centre of~$G$, so that $G/Z=H$. Let $(T,S)$ be a section of $G$ such that $T/S\cong H$. Then $ZT/ZS$ is isomorphic to a quotient of $T/S\cong H$ and is a subquotient of $G/Z=H$, so it is either trivial of the whole of $G/Z$. If $ZT/ZS=\un$, then $ZT=ZS$ and $T=(Z\cap T)S$, so that the group
$$H\cong T/S\cong(Z\cap T)/(Z\cap S)$$
is abelian, contrary to our assumption. Thus $ZT/ZS=G/Z$, hence $ZT=G$ and $S\leq Z$. It follows that $G/S\cong Z/S\times T/S$. Therefore $G$ has an abelian quotient $Z/S$, so $Z/S=\un$ since $G$ is perfect. So we have $S=Z$ and $T=G$. This proves that $(G,Z)$ is the only section in~$\Sigma_H(G)$ and Proposition~\ref{only-section} applies. Therefore
$$S_{H,V}(G) \cong \Tr_1^{G/G}(V)=V \,,$$
as was to be shown.
\endpf

\result{Examples} \label{example} 
{\rm Let $H=A_5$ and let $\varepsilon$ be the sign representation of the group $\Out(A_5)$ of order~2 (over a field $k$ of characteristic different from~2). By Corollaries~\ref{simple} and~\ref{central}, we have $S_{A_5,\varepsilon}(S_5)=0$ and $S_{A_5,\varepsilon}(\tilde A_5)\neq0$. We also have $S_{A_5,\varepsilon}(S_7)=0$, because it is easy to check that $(A_5,\un)$ is the only section in $\Sigma_{A_5}(A_7)$ up to conjugation.}
\fresult

\bigskip
\noindent
Serge Bouc, CNRS-LAMFA, Universit\'e de Picardie - Jules Verne,\\
33, rue St Leu, F-80039 Amiens Cedex~1, France.\\
{\tt serge.bouc@u-picardie.fr}

\medskip
\noindent
Radu Stancu, CNRS-LAMFA, Universit\'e de Picardie - Jules Verne,\\
33, rue St Leu, F-80039 Amiens Cedex~1, France.\\
{\tt radu.stancu@u-picardie.fr}

\medskip
\noindent
Jacques Th\'evenaz, Section de math\'ematiques, EPFL, \\
Station~8, CH-1015 Lausanne, Switzerland.\\
{\tt Jacques.Thevenaz@epfl.ch}


\begin{thebibliography}{1}


\bibitem[BD]{BD}
R.~Boltje, S.~Danz.
\newblock A ghost algebra of the double Burnside algebra in characteristic zero,
\newblock {\em J. Pure Appl. Alg.}, to appear.

\bibitem[Bo1]{Bo1}
S.~Bouc.
\newblock Foncteurs d'ensembles munis d'une double action,
\newblock {\em J. Algebra} 183 (1996), 664--736.

\bibitem[Bo2]{Bo2}
S.~Bouc.
\newblock Biset functors for finite groups,
\newblock Springer Lecture Notes in Mathematics no. 1990 (2010).

\bibitem[BST]{BST}
S.~Bouc, R.~Stancu, J.~Th\'evenaz.
\newblock Simple biset functors and double Burnside ring, 
\newblock {\em J. Pure Appl. Alg.} 217 (2013), 546--566.

\bibitem[BT1]{BT1}
S.~Bouc, J.~Th\'evenaz.
\newblock Gluing torsion endo-permutation modules, 
\newblock {\em J. London Math. Soc.} 78 (2008), 477--501.

\bibitem[BT2]{BT2}
S.~Bouc, J.~Th\'evenaz.
\newblock Stabilizing bisets,
\newblock {Adv. in Math.} 229 (2012), 1610--1639.

\bibitem[TW]{TW}
J.~Th\'evenaz, P.~Webb.
\newblock The structure of Mackey functors,
\newblock {\em Trans. Amer. Math. Soc.} 347 (1995), 1865--1961.

\bibitem[We]{We}
P.~Webb.
\newblock Two classifications of simple Mackey functors with applications to group cohomology and the decomposition of classifying spaces,
\newblock {\em J. Pure Appl. Alg.} 88 (1993), 265--304.


\end{thebibliography}
\end{document}